\documentclass[a4paper,11pt]{article}
\usepackage{amsmath, esint}
\usepackage{amssymb}
\usepackage{amsthm}
\usepackage{bbm}
\usepackage{dsfont}

\usepackage{amsthm}
\usepackage{graphicx}
\usepackage{amssymb}

\usepackage[english]{babel}
\usepackage{latexsym}
\usepackage{mathtools}
\usepackage[arrow, matrix, curve]{xy}
\usepackage{tikz}
\usepackage{float}
\usepackage{enumerate}
\usepackage{leftidx}
\usepackage{longtable}

\newtheorem{thm}{Theorem}
\numberwithin{thm}{section}

\newtheorem{lemma}[thm]{Lemma}

\newtheorem{Prop}[thm]{Proposition}

\newcommand\upperl[5]{\leftidx{^{(#1)\!\!}_{\hphantom{(}#2{\hphantom{)}}\!\!}}{#3}_{#4}^{#5}}

 \numberwithin{equation}{section}

\newcommand{\pibar}{\overline{\Pi}}
\newcommand{\pibarinv}{\overline{\Pi}^{\leftarrow}}

\newcommand{\lambar}{\overline{\Lambda}}

\newcommand{\cF}{\cal F}
\newcommand{\veps}{\varepsilon}

\newcommand{\topr}{\stackrel{\mathrm{P}}{\longrightarrow}}
\newcommand{\todr}{\stackrel{\mathrm{D}}{\longrightarrow}}
\newcommand{\eqdr}{\stackrel{\mathrm{D}}{=}}

\newcommand{\R}{\mathbb{R}}
\newcommand{\D}{\mathbb{D}}
\newcommand{\N}{\mathbb{N}}

\newcommand{\PT}{\mathbb{P}}

\newcommand{\rmd}{{\rm d}}

\newcommand{\halmos}{\quad\hfill\mbox{$\Box$}}

\newcommand{\DD}{\mathbb{D}}

\newcommand{\dto}{\downarrow}

\newcommand{\be}{\begin{equation}}
\newcommand{\ee}{\end{equation}}
\newcommand{\bea}{\begin{eqnarray}}
\newcommand{\eea}{\end{eqnarray}}
\newcommand{\bean}{\begin{eqnarray*}}
\newcommand{\eean}{\end{eqnarray*}}
\newcommand{\ben}{\begin{equation*}}
\newcommand{\een}{\end{equation*}}
\newcommand{\ba}{\begin{aligned}}
\newcommand{\ea}{\end{aligned}}

\newcommand{\bfeins}{{\bf 1}}

\def\nexto{\kern -0.54em}

\def\topr{\buildrel P \over \to }

\newcommand{\PP}{\mathbb{P}}
\newcommand{\EE}{\mathbb{E}}

\newcommand{\PPP}{\mathbb{P}}

\newcommand{\bm}[1]{\mbox{\boldmath $#1$\unboldmath}}

\begin{document}
\title{\bf Small Time Convergence of Subordinators with Regularly or Slowly Varying Canonical Measure}
\author{Ross Maller and Tanja Schindler\thanks{Research
partially supported by ARC Grant DP160103037. \newline
Email: Ross.Maller@anu.edu.au; Tanja.Schindler@anu.edu.au}}


\date{\today}


\maketitle

\begin{abstract}
We consider subordinators  $X_\alpha=(X_\alpha(t))_{t\ge 0}$  in the domain of attraction at 0 of a stable subordinator  $(S_\alpha(t))_{t\ge 0}$ (where $\alpha\in(0,1)$); 
thus, with the property that 
 $\pibar_\alpha$,  the tail function of the canonical measure of $X_\alpha$, 
 is regularly varying of index $-\alpha\in (-1,0)$ as $x\dto 0$.
We also  analyse the boundary case, $\alpha=0$, when  $\pibar_\alpha$ is slowly varying at 0.
When $\alpha\in(0,1)$, we show that  $(t \pibar_\alpha (X_\alpha(t)))^{-1}$ converges in distribution,
as $t\dto 0$, to the random variable $(S_\alpha(1))^\alpha$. 
This latter random variable, as a function of $\alpha$,  converges in distribution as $\alpha\dto 0$ to the inverse of an exponential random variable. 
We prove these convergences, also generalised to functional versions (convergence in $\D[0,1]$),
and to trimmed versions, whereby a fixed number of its largest jumps up to a specified time are subtracted from a process.
The $\alpha=0$ case produces convergence to an extremal process constructed from ordered jumps of a Cauchy subordinator. 
Our results generalise  random walk  and stable process results  of  Darling, Cressie, Kasahara,  Kotani and  Watanabe.
\end{abstract}

\section{Introduction}\label{s1}
A classic result of L\'evy (1937) 
is that stable laws with index $\alpha\in(0,2)$ constitute the entire class of possible non-normal limit laws of a normed and centered random walk in $\R$.
Random walks with such behaviour  are said to be in the domain of attraction of the corresponding stable distribution.

A significant connection, going back to Doeblin (1940), 
and expanded on by Feller (1967, 1971),   
was to  use Karamata's regular variation theory to characterise 
random walks in domains of attraction by regularly varying conditions on the tail of the distribution of the increments of the random walk.
With an appropriate interpretation, the boundary case $\alpha=2$ also corresponds to a stable law, namely the normal distribution, and the corresponding domain of attraction can be characterised
with regular variation-related results. 

What of the other boundary case, $\alpha=0$?
Cressie (1975)   
showed that if $S_\alpha$ is a Stable$(\alpha)$  random variable with index $\alpha\in(0,2)$ and shift constant $\gamma$, then 
$|S_\alpha-\gamma|^\alpha$ converges in distribution as $\alpha \dto 0$ to the reciprocal of an exponential random variable.
Kasahara (1986),   
in a result he attributes to Kotani, extended this in the following way:
let $(S_\alpha(t))_{t\ge 0}$ be a positive stable process of index $\alpha\in(0,1)$, i.e., a subordinator with L\'evy triplet
$(0,0,  x^{-\alpha-1}\rmd x\bfeins_{\{x>0\}})$, having Laplace transform
\ben
\EE e^{-\lambda S_\alpha(t)}= e^{-t\Gamma(1-\alpha)\lambda^\alpha},\ \lambda>0,\ t>0.
\een
Then Kasahara's result states that 
\be\label{Ko}
((S_\alpha(t))^\alpha)_{t\ge 0} \todr (e_t)_{t\ge 0}, \ {\rm as}\ \alpha\dto 0,
\ee
where $\todr$ 
denotes convergence in the Skorohod $J_1$ topology, and $(e_t)$ is an extremal process with marginal distributions
\ben 
\PT\big(e_{t_1}\leq x_1, \ldots, e_{t_n}\leq x_n\big)=
\prod_{i=1}^n  e^{-t_i/x_i},
\een
 for 
 $0\leq t_1<\cdots<t_n$ and $0<x_1<\cdots<x_n$.
 We refer to Resnick (1987)  
for background information on extremal processes. 
 
For each $t>0$, $e_t$ has the distribution of 
 the reciprocal of an exponential  random variable, so \eqref{Ko} represents an extension of the Cressie (1975) result.
The identity
 \be\label{esup}
 e_t \eqdr \sup_{0<s\le t}\Delta \xi_s,
 \ee
 for each $t>0$,   
also holds, where
 $(\xi_t)_{t\ge 0}$ is a Cauchy subordinator, i.e., a  L\'evy process with triplet $(0,0, x^{-2}\rmd x\bfeins_{\{x>0\}})$,
 and jump process $\Delta \xi_t:=\xi_t-\xi_{t-}$, $t>0$.
 
 When $0<\alpha<2$, the tail of the increment distribution of a random walk in  the domain of attraction of a Stable$(\alpha)$ distribution is  regularly varying  at $\infty$ with index $-\alpha$.
So for  the  boundary case, $\alpha=0$, it is natural to consider a slowly varying tail. In this case affine norming and centering of the random walk cannot lead to a finite nondegenerate limit  random variable, but a transformation,  whereby the tail of the increment distribution is applied as a function to the random walk, and norming is by the sample size, produces
as a limiting  random variable the reciprocal of an exponential  random variable.
This was proved by Darling (1952)  
in a 1-dimensional version, and, subsequent to this, in Watanabe (1980),  
the random walk is interpolated to a function in $\D[0,1]$, and finite dimensional convergence of the resulting process is proved. 
In a later paper,  Kasahara (1986)  
proved $J_1$ convergence of the functions in $\D[0,1]$.

In view of this background,  the continuous time environment is a natural one in which to consider results like these, and the aim of the present paper is, firstly to transfer from random walk versions to L\'evy processes, in which the convergence is for small time parameter, rather than large time, and, secondly, to generalise the results to trimmed versions of L\'evy processes.
By ``trimming" we mean removing a fixed number of large jumps of the processes.
This is natural in the random walk context, because the slowly varying, heavy tails are associated with large jumps (``outliers") in the random walk, and it is 
interesting in the process context as the effect of a slowly varying measure  near 0 is previously little explored. 
 Apart from these aspects, some quite interesting analytical differences occur between the small and large time situations.

Thus our basic assumption will be of the kind that
a generic L\'evy process $(Y(t))_{t\ge 0}$ 
with triplet $(\gamma_Y, \sigma^2_Y, \Pi_Y(\rmd y))$, is in a non-normal domain of attraction
at small times, by which we mean there exist non-stochastic functions $a_t \in \R$ and $b_t > 0$ such that 
\be\label{dom}
 \frac{Y(t) - a_t}{b_t} \todr S,\ {\rm as }\ t \dto 0,
\ee 
where $S$ is an almost surely (a.s.) finite, non-degenerate, non-normal  random variable.

Conditions on the L\'evy measure for \eqref{dom} to hold (in small time) can be deduced from Theorem 2.3 of Maller and Mason (2008),  
whose result can also be used to show that \eqref{dom} can be extended to convergence in $\DD[0,1]$;
that is, 
\be\label{funcFC1}
\left( \frac{Y(\lambda t) - \lambda a_t}{b_t} \right)_{0<\lambda \le 1}
\to (S(\lambda))_{0<\lambda \le 1},\ {\rm as }\ t \dto 0,
\ee
weakly with respect to the Skorohod $J_1$ topology.
Then \eqref{dom} is equivalent to the two-sided tail $\pibar_Y$ of $Y$ being regularly varying at 0   with index $\alpha\in(0,2)$, together with a balance condition on the right and left tails of the L\'evy measure $\Pi_Y$.
The limit  random variable $S$ in \eqref{dom} has the distribution of  $S_\alpha(1)$,  where $(S_\alpha(\lambda)_{0<\lambda \le 1}$
is a Stable($\alpha$) L\'evy process.

In Buchmann, Ipsen and Maller (2017)  
\eqref{funcFC1} was extended to a functional theorem for a trimmed version of $Y$, which result will be quoted below (see the proof of Theorem \ref{th2}).
The case of a slowly varying tail for $\Pi$ seems not to have been considered before, in our context  (but see Kevei and Mason (2014) 
and Ipsen, Maller and Resnick (2018)  
for limits of ratios of large jumps of subordinators in this case).
Although stated in \eqref{dom} and \eqref{funcFC1} for general L\'evy processes, from now on we restrict ourselves to subordinators. Some discussion relevant to this is given at the end of the next section.


\section{Notation and Statement of Results}\label{s2} 
All processes will be defined on a probability space $(\Omega, \cF, \PPP)$.
Since the index $\alpha$ will be a variable in  our results, we  have to indicate its presence in the notation.
We have tried to come up with a notation that is minimal but clear and conveys the necessary information.

For each $\alpha\ge 0$ let $\left(X_\alpha(t)\right)_{t\ge 0}$ be a driftless  subordinator with canonical measure
$\Pi_\alpha(\rmd x)$, where $\Pi_\alpha$ 
has tail $\pibar_\alpha(x):=
\Pi_\alpha\{(x,\infty)\}$, $x>0$,  satisfying
$\pibar_\alpha(x)=x^{-\alpha}L(x){\bf 1}_{\{x>0\}}$, with $0\le \alpha<1$, and $L(x)$  a function slowly varying as $x\dto 0$.
For the $\alpha=0$ case, simply write  $X(t):=X_0(t)$ and $\Pi:=\Pi_0$. In this case, $L(x)$ is assumed to be nonincreasing with $L(0+)=\infty$.
Since the  processes $X_\alpha(t)$ are subordinators, $\alpha$ is necessarily restricted to  $[0,1)$.

Our development goes as follows. For each $\alpha\in(0,1)$, $X_\alpha(t)$ is
in the domain of attraction of a positive Stable$(\alpha)$ distribution as $t\dto 0$;
 in fact,  the process 
$(X_\alpha(\lambda t))_{0<\lambda\le 1}$,  converges in $\D[0,1]$, as $t\dto 0$, after norming, 
to  a  Stable$(\alpha)$ process  $(S_\alpha(\lambda))_{\lambda\ge 0}$.
This implies that
$\big(t\pibar_\alpha
\big(X_\alpha(t\lambda)\big)\big)^{-1}$ converges to 
$(S_\alpha(\lambda)\big)^{\alpha}$ in $\D[0,1]$.
In turn, this latter process
itself converges in distribution, as $\alpha\dto 0$,
to the largest jump up till time $\lambda$ of a  Stable$(1)$ (Cauchy) process with measure $x^{-2}\rmd x{\bf 1}_{\{x>0\}}$.
We denote this process as $(\xi_t)_{t\ge 0}$, consistent with the notation in \eqref{esup}.

These results are included in our main theorem, 
Theorem \ref{th1}, set out in diagrammatic form below.
It deals, not just with the processes mentioned, but also with ``trimmed" versions of them.
To introduce trimmed processes,
 write $(\Delta X_\alpha(t):=X_\alpha(t)-X_\alpha(t-))_{t>0}$, with   $\Delta X_\alpha(0)=0$, for the jump process of $X_\alpha$, and
 $\Delta X_\alpha^{(1)}(t) \ge \Delta X_\alpha^{(2)}(t)\ge \cdots$
 for the ordered jumps at time $t>0$.
 Since $\Pi\{(0,\infty)\}=\infty$,
there are infinitely many positive jumps, a.s., in any finite time interval $[0,t]$, $t>0$, 
the $\Delta X_\alpha^{(r)}(t)$ are positive a.s.\ for all $t>0$, and 
$\lim_{t\dto 0} \Delta X_\alpha^{(r)}(t)=0$ a.s.\ for all $r\in\N$.
(Throughout, let $\N:=\{1,2,\ldots\}$ and $\N_0:=\{0,1,2,\ldots\}$.)
 The {\it $r$-trimmed process} is defined to be $X_\alpha(t)$ minus its $r$  largest jumps, at a given time $t$:
\begin{equation}\label{trims}
{}^{(r)}X_\alpha(t):= X_\alpha(t)-\sum_{i=1}^r {\Delta X}_\alpha^{(i)}(t), \ r\in\N, \ t>0
\end{equation}
(and we set $^{(0)}X_\alpha(t)\equiv X_\alpha(t)$).
Detailed definitions and properties of this kind of ordering  and trimming are given in Buchmann, Ipsen and Maller (2016), 
where the (positive) $\Delta X_\alpha(t)$ are identified with the points of a Poisson point process 
on $[0,\infty)$.

We similarly denote the 
ordered jumps up till time $\lambda$ 
of the Cauchy process  $(\xi_\lambda)_{\lambda\ge 0}$ with
 jump process $(\Delta \xi_\lambda)_{\lambda\ge 0}$
as $\Delta \xi_\lambda^{(1)}\ge \Delta \xi_\lambda^{(2)}\ge \cdots$.

\begin{thm}\label{th1}
For each  $\alpha\in[0,1)$ let $(X_\alpha(t))$ be a driftless subordinator whose tail measure 
$\pibar_\alpha$ is regularly varying at zero with exponent $-\alpha$ and satisfies $\pibar_\alpha(0+)=\infty$; and for each $r\in\N_0$ let $({}^{(r)}X_\alpha(t))$ be the trimmed version of $(X_\alpha(t))$ defined in \eqref{trims}. 
When $\alpha=0$ assume in addition 
\begin{align}\label{eq: Xt Delta Xt as}
 \lim_{t\to 0} \frac{^{(r)} X_0(t)} {\Delta X_0^{(r+1)}(t)}=1,\ {\rm  a.s.,}
\end{align}
where ``a.s." denotes almost sure convergence.
Then for all $r\in\mathbb{N}_0$ we have the following convergences in distribution, as $t\dto 0$, and as $\alpha\dto 0$,  with respect to the Skorohod $J_1$-topology and the parameter $\lambda\in(0,1]$:
\begin{figure}[H]
\centering
\begin{tikzpicture}
\node  (ul)  at (0,4)  {$\big(t\pibar_\alpha
\big(^{(r)}X_\alpha(t\lambda)\big)\big)^{-1}$};
\node  (ur)  at (5,4)  {$\big(t\overline{\Pi}
\big(^{(r)}X_0(t\lambda)\big)\big)^{-1}$};
\node  (ll)  at (0,1)  {$\big({}^{(r)}S_\alpha(\lambda)\big)^{\alpha}$};
\node  (lr)  at (5,1)  {$(\Delta\xi_{\lambda}^{\left(r+1\right)})$};

\draw[->] 
(ul) edge node [midway, sloped, above] {$t\to 0$}
         node [midway,  below] {\,\,\,\,\,\,\,}(ll)
(ur) edge node [midway, sloped, above] {$t\to 0$} 
          node [midway,  below] {\,\,\,\,\,\,\,}(lr)
(ll) edge node [midway, sloped, above] {$\alpha\to 0$} 
          node [midway, sloped, below] {}(lr);
\end{tikzpicture}
\caption{{\it Main Convergence Diagram.
The upper nodes  represent  processes in $0<\lambda\le 1$, indexed by $t>0$.
The lower nodes  represent  processes in $0<\lambda\le 1$. 
The index $r\in\N_0$  indicates the order of trimming. 
The vertical arrows indicate process convergence of the upper node processes as $t\dto 0$ to the lower node processes
for each $\alpha\in(0,1)$ on the left, and with $\alpha$ set equal to 0 on the right.
The horizontal arrow indicates process convergence of the left lower node process as $\alpha\dto 0$ to the right lower node process.}
}\label{diagram}
\end{figure}
\end{thm}

\medskip\noindent{\bf Remarks.}\ (i)\ 
Some comment on Figure \ref{diagram} is in order.
Since $\pibar_\alpha(0+)=\infty$ (i.e., $\Pi$ is of ``infinite activity") for each $\alpha\ge 0$, 
and $\lim_{t\dto 0} {}^{(r)}X_\alpha(t)=0$ a.s.,
we have
$\lim_{t\dto 0} \pibar_\alpha( {}^{(r)}X_\alpha(t))=\infty$ a.s., and under the regularly varying (at 0) assumption 
we impose on $\pibar_\alpha$, it turns out that multiplying by $t$ is the correct scaling to get a 
nondegenerate limit law for $ \pibar_\alpha( {}^{(r)}X_\alpha(t))$ as $t\dto 0$. 
It is then convenient to consider the limit of the reciprocal of $ t\pibar_\alpha( {}^{(r)}X_\alpha(t))$ 
as we do in the  topmost  entries of Figure \ref{diagram} because it produces the trimmed stable in the upright orientation as we see in the bottom left entry of the figure, thereby  providing a direct generalisation of the Kotani result in \eqref{Ko}.
Taking the function $\pibar_\alpha( {}^{(r)}X_\alpha(t))$ of $ {}^{(r)}X_\alpha(t)$ is a natural way of 
generalising the Darling (1952) result for random walks, but it's clear that some quite different considerations enter in;
note for example that $\pibar$ slowly varying at zero reflects a mild singularity, while $\alpha\in(0,1)$ is steeper --
whereas, at infinity, a slowly varying $\pibar$ betokens a very heavy tailed random walk. 

(ii)\ 
The appearance of the almost sure condition \eqref{eq: Xt Delta Xt as} among the other weak convergence results is at first
surprising. We discuss this in more detail after the proof of Theorem \ref{th4.1}. 

(iii) Given the exposition in \eqref{dom} and \eqref{funcFC1}, it is logical to ask if there are versions of the convergences in Theorem \ref{th1} for (necessarily centered) general L\'evy processes, other than subordinators. We have not investigated in detail whether this can be done, 
but the results for subordinators are certainly of interest in themselves, (i) as being generalisations of non-negative random walk versions which have appeared in the literature discussed in Section \ref{s1}, and, (ii)  because subordinators 
and their jumps play a prominent role for example in  the theory of Poisson-Dirichlet distributions initiated by Kingman (1975), 
which is not geared to the application of general L\'evy processes. 
A further interesting point is that the  Kingman Poisson-Dirichlet development relates at its heart to the {\it small time} behaviour of the stable subordinators, such as we consider here.

\section{Convergence of $X_\alpha(t)$ as $t\dto 0$, for fixed $\alpha\in(0,1)$}\label{sec: ad a} 
In this section we prove the lefthand vertical convergence in Figure \ref{diagram}.
Here the parameter 
$\alpha$ does not vary; the convergence is as $t\dto 0$, for fixed $\alpha$. 

\begin{thm}\label{th2}
Fix  $\alpha\in(0,1)$ and let $(X_\alpha(t))$ be a driftless subordinator whose tail measure 
$\pibar_\alpha$ is regularly varying at zero with exponent $-\alpha$. For each $r\in\N$ let $({}^{(r)}X_\alpha(t))$ be the trimmed process defined in \eqref{trims}. 
Then
\begin{align}\label{5}
\left(
\frac{1}
{t\overline{\Pi}_\alpha\left(^{(r)}X_\alpha(t\lambda)\right)}
\right)_{0<\lambda\leq 1}
\todr
\big(\big(^{(r)}S_\alpha(\lambda)\big)^\alpha\big)_{0<\lambda\leq 1},\ {\rm as}\ t \dto 0,
 \end{align}
  with respect to the $J_1$-topology.
\end{thm}

In what follows, define  the generalized inverse function of a monotonically decreasing function $g$ by
$g^{\leftarrow}\left(x\right)\coloneqq \inf\left\{y>0\colon g\left(y\right)\leq x\right\}$,  for $x>0$.

\bigskip\noindent{\bf Proof of Theorem \ref{th2}:}\
Fix  $\alpha\in(0,1)$ and assume 
$\pibar_\alpha$ is regularly varying at zero with exponent $-\alpha$.
Then 
\begin{align} \label{eq: conv in J1}
 \left(\frac{X_\alpha(t\lambda)}
 {\pibarinv_\alpha \left(1/t\right)}\right)_{0<\lambda\leq 1}\overset{D}{\longrightarrow}
 \left(S_\alpha(\lambda)\right)_{0<\lambda\leq 1},\ {\rm as}\ t \dto 0,
 \end{align}
  with respect to the $J_1$-topology.
  This result follows 
   from Theorem 2.3 of Maller and Mason (2008)   
(see also their references for antecedents) quoted as the generic version in \eqref{funcFC1}.   
  Maller and Mason  do not mention that the norming function $b_t$ in \eqref{funcFC1} can be taken as the inverse function to the tail measure of the process, or that, in the driftless subordinator case, the centering function $a_t$ can be taken as 0, as we have done in \eqref{eq: conv in J1}; but these facts are easily checked.
  
  Taking \eqref{eq: conv in J1} as given, it 
  further implies the trimmed version 
\begin{align}\label{eq: Xalpha Pi t}
 \left(\frac{^{(r)}X_\alpha(t\lambda)}{\pibar_\alpha^{\leftarrow}\left(1/t\right)}\right)_{0<\lambda\leq 1}\overset{D}{\longrightarrow}\big(^{(r)}S_\alpha(\lambda)\big)_{0<\lambda\leq 1},\ {\rm as}\ t \dto 0,\ {\rm for}\ r\in\N,
\end{align}
 with respect to the $J_1$-topology, as shown in Theorem 3 of Buchmann, Ipsen and  Maller (2017).
The convergence in \eqref{eq: Xalpha Pi t} additionally  implies  
\begin{align}\label{10}
 \left(
\frac{\pibar_\alpha\left(\pibarinv_\alpha\left(1/t\right)\right)}
{\overline{\Pi}_\alpha\left(^{(r)}X_\alpha(t\lambda)\right)}
\right)_{0<\lambda\leq 1}
\todr
\left(\big(^{(r)}S_\alpha(\lambda)\big)^\alpha\right)_{0<\lambda\leq 1},
\end{align}
 by application of the following Lemma \ref{lem1}, 
and \eqref{10} implies \eqref{5}, thereby completing
the  proof of Theorem \ref{th2}.  \halmos

\begin{lemma}\label{lem1}
Suppose 
$\pibar_\alpha$ is regularly varying at zero with exponent $-\alpha$, $\alpha>0$. Then
for two functions $f_t>0$ and $g_t>0$ on $[0,\infty)$ with
$\lim_{t\dto0}f_t=\lim_{t\dto0}g_t=0$,
we have 
 $\lim_{t\dto 0}\pibar_\alpha\left(g_t\right)/\pibar_\alpha\left(f_t\right)=c^{\alpha}$
 if and only if
 $\lim_{t\dto0}f_t/g_t=c\in (0,\infty)$.
\end{lemma}

\bigskip\noindent{\bf Proof of Lemma \ref{lem1}:}\
This is a straightforward application of  Potter's bounds, see for example Theorem 1.5.6 of Bingham, Goldie and Teugels (1987).
We omit the details.
\halmos

\section{Convergence of $X_t=X_0(t)$ as $t\dto 0$, Case $\alpha= 0$}\label{sec: ad b}
Next we prove the righthand vertical convergence in Figure 1. The
process $X(t)=X_0(t)$  is now assumed to have tail $\pibar(x)$  slowly varying as $x\dto 0$, and the results in this section formally correspond to the case $\alpha=0$.
So we drop the subscript $\alpha$ and write $X_t$ rather than $X(t)$
throughout this section.  Keep $r\in\N_0$ fixed.
Recall that   $\Delta \xi_\lambda^{(1)}\ge \Delta \xi_\lambda^{(2)}\ge \cdots$
are the  ordered jumps, up till time $\lambda$, of $\xi_\lambda$. 
The main result for  this section is:

\begin{thm}\label{th4.1}
Suppose $X_t$ is a driftless subordinator whose L\'{e}vy measure $\Pi$ has tail $\overline{\Pi}$  slowly varying at zero. 
Assume \eqref{eq: Xt Delta Xt as} in addition.
Then 
\begin{align*}
\left(
\frac{1}{t \pibar(^{(r)}X_{t\lambda})}
\right)_{0<\lambda\leq 1}
\overset{D}{\longrightarrow}  \big(\Delta\xi_{\lambda}^{(r+1)}\big)_{0<\lambda\leq 1},\ {\rm as}\ t \dto 0,
 \end{align*}
  with respect to the $J_1$-topology. 
\end{thm}
Proof of Theorem \ref{th4.1} proceeds by way of some lemmas and propositions. 
The first lemma proves convergence in the supremum norm of the difference of two quantities to 0, stronger than proving $J_1$ convergence.

\begin{lemma}\label{lem2} 
Assume the conditions of Theorem \ref{th4.1}, including \eqref{eq: Xt Delta Xt as}.
Then for each $r\in\N_0$ 
\begin{align}\label{eq: to prove}
\sup_{0<\lambda\leq 1}
\left|\frac{1}{t \overline{\Pi}\left(^{(r)}X_{t\lambda}\right)}
-\frac{1}{t \overline{\Pi}\big(\Delta X_{t\lambda}^{\left(r+1\right)}\big)}\right|\topr 0, \ {\rm as}\ t\dto 0.
\end{align}
\end{lemma}

\bigskip\noindent{\bf Proof of Lemma \ref{lem2}:}\
Hold  $r\in\N_0$ fixed throughout. 
Since $X_t$ is a subordinator, its jumps are positive, and so 
$\upperl{r}{}{X}{t\lambda}{}\geq \Delta X_{t\lambda}^{\left(r+1\right)}$ for $t>0$ and $\lambda\in(0,1]$.
Thus $1/\overline{\Pi}(\upperl{r}{}{X}{t\lambda}{})\geq
1/\overline{\Pi}\big(\Delta X_{t\lambda}^{\left(r+1\right)}\big)$, and for
\eqref{eq: to prove} it suffices to prove that
 for all $y>0$ and $ \eta>0$ there exists $t_0=t_0(y,\eta)>0$ such that $t\in(0,t_0)$ implies
\begin{align}\label{16}
 \mathbb{P}\left(\sup_{0<\lambda\leq 1}
 \left(\frac{1}{t\overline{\Pi}(\upperl{r}{}{X}{t\lambda}{})}
-\frac{1}{t \overline{\Pi}\big(\Delta X_{t\lambda}^{(r+1)}\big)}\right)>y\right)<\eta.
\end{align}
Take  $K>0$. The  left hand side  of \eqref{16} equals 
\ben
 \mathbb{P}\left(\sup_{0<\lambda\leq 1}\frac{1}
{t \overline{\Pi}\big(\Delta X_{t\lambda}^{\left(r+1\right)}\big)} 
\left(\frac{t \overline{\Pi}\big(\Delta X_{t\lambda}^{\left(r+1\right)}\big)}
{t \overline{\Pi}\big(\upperl{r}{}{X}{t\lambda}{}\big)}-1\right)>y\right),
\een
and this is bounded above by
\bea\label{eq: t Pi X - t Pi Delta X}
&&
\mathbb{P}\left(\sup_{0<\lambda\leq 1}\frac{1}{t \overline{\Pi}\big(\Delta X_{t\lambda}^{\left(r+1\right)}\big)}\leq K,\,
\sup_{0<\lambda\leq 1}\left(\frac{ \overline{\Pi}\big(\Delta X_{t\lambda}^{\left(r+1\right)}\big)}{ \overline{\Pi}(\upperl{r}{}{X}{t\lambda}{})}-1\right)>
\frac{y}{K}\right)\cr
&&
\hskip4cm 
+\mathbb{P}\left(\sup_{0<\lambda\leq 1}\frac{1}{t \overline{\Pi}\big(\Delta X_{t\lambda}^{\left(r+1\right)}\big)}>K\right).
\eea

We bound  the first probability in
\eqref{eq: t Pi X - t Pi Delta X}  by ignoring the first supremum in it. To deal with the remaining part of that term, we need to invoke \eqref{eq: Xt Delta Xt as}.
This condition implies 
that there is an event $\Omega_1\subseteq \Omega$ with $\PP(\Omega_1)=1$ such that, for 
 $\omega\in\Omega_1$ and $\delta>0$, there exists $t_1\in(0,t_0)$ such that for $t\in(0,t_1)$ we have $W_t^{(r)}\coloneqq \upperl{r}{}{X}{t}{}/\Delta X_{t}^{\left(r+1\right)}<1+\delta$, and thus $\sup_{0<\lambda\leq 1} W_{t\lambda}^{(r)}<1+\delta$.
Hence, we can find $t_2\in(0,t_1)$ such that, for  $t\in(0,t_2)$,
\begin{align*}
 \mathbb{P}\Big(\sup_{0<\lambda\le 1}W_{t\lambda}^{(r)}>2\Big)<\frac{\eta}{3}.
\end{align*}
Then for $t\in(0,t_2)$ we have  
\begin{align}\label{eq: t Pi X - t Pi Delta X 2a}
&
\mathbb{P}\left(\sup_{0<\lambda\le 1}\left(\frac{\overline{\Pi}\big(\Delta X_{t\lambda}^{\left(r+1\right)}\big)}{\overline{\Pi}(\upperl{r}{}{X}{t\lambda}{})}-1\right)>\frac{y}{K}\right)\notag\\
&\leq
 \mathbb{P}\left(\sup_{0<\lambda\le 1}\left(\frac{\overline{\Pi}\big(\Delta X_{t\lambda}^{\left(r+1\right)}\big)}{\overline{\Pi}\big(2\Delta X_{t\lambda}^{\left(r+1\right)}\big)}-1\right)>\frac{y}{K},\ \sup_{0<\lambda\le 1}W_{t\lambda}^{(r)}\leq 2\right)\cr
&\hskip7cm 
+
\mathbb{P}\Big(\sup_{0<\lambda\le 1}W_{t\lambda}^{(r)}>2\Big)\notag\\
&\leq
 \mathbb{P}\left(\sup_{0<\lambda\le 1}\left(\frac{\overline{\Pi}\big(\Delta X_{t\lambda}^{\left(r+1\right)}\big)}
{\overline{\Pi}\big(2\Delta X_{t\lambda}^{\left(r+1\right)}\big)}-1\right)>\frac{y}{K}\right)+\frac{\eta}{3}.
\end{align}
The slow variation of $\overline{\Pi}$  implies 
there exists $x_0>0$ such that, for all $x\in(0,x_0]$,
$\overline{\Pi}\left(x\right)/\overline{\Pi}\left(2x\right)-1\leq y/K$.
Further, notice that $\{\Delta X_{t}^{\left(r+1\right)}\leq x_0\}$ implies
$\{\sup_{0<\lambda\le 1}\Delta X_{t\lambda}^{\left(r+1\right)}\leq x_0\}$,
and thus, when $\Delta X_{t}^{\left(r+1\right)}\leq x_0$,
\begin{align*}
\left\{\sup_{0<\lambda\le 1}\left(\frac{\overline{\Pi}\big(\Delta X_{t\lambda}^{\left(r+1\right)}\big)}
{\overline{\Pi}\big(2\Delta X_{t\lambda}^{\left(r+1\right)}\big)}-1\right)\leq \frac{y}{K}\right\}.
\end{align*}
Hence, the probability on the righthand side of \eqref{eq: t Pi X - t Pi Delta X 2a} can be estimated as
\begin{align}\label{eq: t Pi X - t Pi Delta X 2a1}
 \mathbb{P}\left(\sup_{0<\lambda\le 1}\left(\frac{\overline{\Pi}\big(\Delta X_{t\lambda}^{\left(r+1\right)}\big)}
{\overline{\Pi}\big(2\Delta X_{t\lambda}^{\left(r+1\right)}\big)}-1\right)>\frac{y}{K}\right)
&\leq 
\mathbb{P}\Big(\sup_{0<\lambda\leq 1}\Delta X_{t\lambda}^{\left(r+1\right)}> x_0\Big)\notag\\
&=\mathbb{P}\big(\Delta X_{t}^{\left(r+1\right)}> x_0\big).
\end{align}
Since $\lim_{t\dto 0} \Delta X_t^{(r)}=0$ a.s.,
there exists $t_3\in(0,t_2)$ such that the righthand side of 
\eqref{eq: t Pi X - t Pi Delta X 2a1} does not exceed $\eta/3$, for $t\in(0, t_3)$.

To estimate the second probability on the  righthand side  of \eqref{eq: t Pi X - t Pi Delta X},  we will use that
there exists $K>0$ and $t_4\in(0,t_3)$ such that, for $t\in(0, t_4)$,
\begin{align}\label{eq: t Pi X - t Pi Delta X 1}
\mathbb{P}\left(\sup_{0<\lambda\leq 1}
\frac{1}{t \overline{\Pi}\big(\Delta X_{t\lambda}^{\left(r+1\right)}\big)}>K\right)
\leq \mathbb{P}\left(\frac{1}{t \overline{\Pi}\big(\Delta X_{t}^{\left(r+1\right)}\big)}>K\right)
\leq \frac{\eta}{3}.
\end{align}
This holds because,  as a special case of the convergence in Proposition \ref{prop4.3} below, 
$ 1/ t \overline{\Pi}\big(\Delta X_{t}^{\left(r+1\right)}\big)$ converges to a finite positive random variable; we defer proof of \eqref{eq: t Pi X - t Pi Delta X 1}
till then.

Accepting \eqref{eq: t Pi X - t Pi Delta X 1}, then, we can combine
 \eqref{eq: t Pi X - t Pi Delta X} with \eqref{eq: t Pi X - t Pi Delta X 2a}, 
 \eqref{eq: t Pi X - t Pi Delta X 2a1} and 
 \eqref{eq: t Pi X - t Pi Delta X 1} to get, for $t\in(0,t_4)$,
\begin{align*}
\mathbb{P}\left(\sup_{0<\lambda\leq 1}\left(\frac{1}{t\overline{\Pi}\left(^{(r)}X_{t\lambda}\right)}
-\frac{1}{t \overline{\Pi}\big(\Delta X_{t\lambda}^{\left(r+1\right)}\big)}\right)>y\right)\leq 3\Big( \frac{\eta}{3}\Big)=\eta.
\end{align*}
Since $\eta$ is arbitrary this completes the proof of \eqref{16}, and  of Lemma \ref{lem2}.
\halmos

Now write
\begin{align*}
\frac{1}{t\overline{\Pi}\left(^{(r)}X_{t\lambda}\right)}=\left(\frac{1}{t\overline{\Pi}\left(^{(r)}X_{t\lambda}\right)}
-\frac{1}{t \overline{\Pi}\big(\Delta X_{t\lambda}^{\left(r+1\right)}\big)}\right)
+\frac{1}{t \overline{\Pi}\big(\Delta X_{t\lambda}^{\left(r+1\right)}\big)}. 
\end{align*}
By  Lemma \ref{lem2} the first summand converges to zero in probability  uniformly in $0<\lambda\le 1$. 
Thus, the processes
\begin{align*}
\left(
\frac{1}{t\overline{\Pi}\left(^{(r)}X_{t\lambda}\right)}
\right)_{0<\lambda\le 1}
\quad {\rm and}\quad 
\left(
\frac{1}{t \overline{\Pi}\big(\Delta X_{t\lambda}^{\left(r+1\right)}\big)}
\right)_{0<\lambda\le 1}
\end{align*}
have the same limit in distribution as $t\dto 0$.
So to complete the proof of Theorem \ref{th4.1} it remains only to prove the following proposition.

\begin{Prop}\label{prop4.3}
Assume the conditions of Theorem \ref{th4.1}, including \eqref{eq: Xt Delta Xt as}. Then, for all $r\in\mathbb{N}$, as $t\dto 0$, 
\begin{align}\label{21}
\left(\frac{1}{t\overline{\Pi}\big(\Delta X_{t\lambda}^{\left(r\right)}\big)}\right)_{0<\lambda\leq 1}\overset{D}{\longrightarrow}
\big(\Delta\xi_{\lambda}^{\left(r\right)}\big)_{0<\lambda\leq 1}, \ {\rm in}\ \D[0,1].
\end{align}
\end{Prop}
We prove this in a classical way, first establishing
 finite dimensional (``fidi") convergence, 
then tightness of the process on the left of \eqref{21}.
This is done in the next two subsections.

\subsection{Proof of  fidi convergence  in Proposition \ref{prop4.3}}\label{sec: fin dim conv}
Define  the following random variables
\be\label{Zdef}
Z_{r,t,\lambda}\coloneqq \frac{1}{t \overline{\Pi}
\big(\Delta X^{(r)}_{t\lambda}\big)}, \ r\in\mathbb{N}, \  t>0, \ \lambda>0,
\ee
and note that $Z_{r,t,\lambda}$ is nondecreasing in $\lambda$. 
Recall that  $\Delta \xi_\lambda^{(1)}\ge \Delta \xi_\lambda^{(2)}\ge \cdots$ are the  ordered jumps, at time $\lambda$, of the Cauchy process  $(\xi_\lambda)_{\lambda\ge 0}$ having  L\'evy  measure $x^{-2}\rmd x\bfeins_{\{x>0\}}$.
Let $\lambda_1<\cdots<\lambda_n$.
We aim to show
\begin{align}\label{eq: fin dim distr to show}
 \MoveEqLeft\lim_{t\dto 0}\mathbb{P}\left(Z_{r,t,\lambda_1}\leq y_1,\ldots, Z_{r,t,\lambda_n}\leq y_n\right)\notag\\
 &=
 \mathbb{P}\big(\Delta\xi^{\left(r\right)}_{\lambda_1}\leq y_1,\ldots, \Delta\xi^{(r)}_{\lambda_n}\leq y_n\big),\
n, r\in\N, 
\end{align}
wherein it is sufficient to restrict ourselves to values $0<y_1<\cdots< y_n$, since $\left\{ Z_{r,t,\lambda_i}\leq y_i\right\}\supseteq\left\{ Z_{r,t,\lambda_j}\leq y_j\right\}$ whenever $i<j$ and $y_i\geq y_j$. 

For formal reasons let $\lambda_0\coloneqq 0$ and  $y_{n+1}\coloneqq\infty$,  and introduce  triangular arrays of random variables $\left(V_{\ell,j}\right)_{1\leq \ell\leq j\leq n}$  and $(\widetilde{V}_{\ell,j,t})_{1\leq \ell\leq  j\leq n,t\geq 0}$ by  setting 
\be\label{eq: def V ij}
V_{\ell,j}
\coloneqq\#\left\{s\in\left(\lambda_{\ell-1},\lambda_{\ell}\right]\colon \Delta \xi_s\in\left(y_j,y_{j+1}\right]\right\}
\ee
and 
\ben 
\widetilde{V}_{\ell,j,t}
\coloneqq\#\left\{s\in \left(t\lambda_{\ell-1},t\lambda_{\ell}\right]\colon \Delta X_s\in\left(\overline{\Pi}^{\leftarrow}
\left((t y_{j})^{-1}\right),\,\overline{\Pi}^{\leftarrow}\left(t y_{j+1})^{-1}\right)\right]\right\},
\een
for $t>0$ and pairs $\ell, j$ fulfilling $1\leq \ell\leq j\leq n$.
The events
 $\{\Delta\xi_{\lambda_{i}}^{\left(r\right)}\leq y_{i}\}$
and 
$\{\sum_{\ell=1}^{i}\sum_{j=i}^{n}V_{\ell,j}\leq r-1\}$
 are equal. This can be seen as follows. By the definition of $V_{\ell,j}$ we have that
 $\sum_{j=i}^{n}V_{\ell,j}=\#\left\{s\in\left(\lambda_{\ell-1},\lambda_{\ell}\right]\colon \Delta \xi_s>y_i\right\}$. 
Thus,
\begin{align*}
 \sum_{\ell=1}^{i}\sum_{j=i}^{n}V_{\ell,j}
 &=\sum_{\ell=1}^{i}\#\left\{s\in\left(\lambda_{\ell-1},\lambda_{\ell}\right]\colon \Delta \xi_s>y_i\right\}\\
 &=\#\left\{s\in\left(0,\lambda_i\right]\colon \Delta \xi_s>y_i\right\}.
\end{align*}
Hence, $\sum_{\ell=1}^{i}\sum_{j=i}^{n}V_{\ell,j}\leq r-1$ holds if and only if
$\#\left\{s\in\left(0,\lambda_i\right]\colon \Delta \xi_s>y_i\right\}\leq r-1$, 
which is equivalent to $\{\Delta\xi_{\lambda_{i}}^{\left(r\right)}\leq y_{i}\}$.

We assert that the event on the right hand side of \eqref{eq: fin dim distr to show}
can be written as a finite union of disjoint events, each of which is the intersection of a finite number of events 
 of the form $\left\{V_{\ell,j}=\kappa_{\ell,j}\right\}$. 
Here the  $\left(\kappa_{\ell,j}\right)_{1\leq \ell\leq j\leq n}$ are triangular arrays of  non-negative integers in which  the $V_{\ell,j}$ and $\widetilde{V}_{\ell,j,t}$ take values.
To verify that assertion,
define
\begin{align*}
 B_{r,n,i}\coloneqq\bigg\{\bm{\kappa}=(\kappa_{\ell,j})_{1\leq \ell\leq j\leq n}\colon \sum_{\ell=1}^{i}\sum_{j=i}^{n}\kappa_{\ell,j}\leq r-1\bigg\}.
\end{align*} 
Assume that for a given tuple $\bm{\kappa}=(\kappa_{\ell,j})$ we have that
$\left\{V_{\ell,j}=\kappa_{\ell,j}\right\}$ for all pairs $\ell, j$
with
$1\leq \ell\leq j\leq n$.
Then $\sum_{\ell=1}^{i}\sum_{j=i}^{n}V_{\ell,j}\leq r-1$ holds if and only if 
$\bm{\kappa}\in B_{r,n,i}$.
On the other hand,  that the event $\left\{V_{\ell,j}=\kappa_{\ell,j}\right\}$ holds simultaneously 
for all pairs $\ell, j$ with
$1\leq \ell\leq j\leq n$,  can also be written as 
$\bigcap_{1\leq \ell\leq j\leq n}\left\{V_{\ell,j}=\kappa_{\ell,j}\right\}$.
This implies 
\begin{align*}
\big\{\Delta\xi_{\lambda_{i}}^{\left(r\right)}\leq y_{i}\big\}
&=
\bigg\{\sum_{\ell=1}^{i}\sum_{j=i}^{n}V_{\ell,j}\leq r-1\bigg\} \cr
&=
\bigcup_{\bm{\kappa}=(\kappa_{\ell,j})\in B_{r,n,i}} \bigcap_{1\leq \ell\leq j\leq n}\left\{V_{\ell,j}=\kappa_{\ell,j}\right\}.
\end{align*}
Now let $A_{r,n}\coloneqq \bigcap_{i=1}^n B_{r,n,i}$, so that $A_{r,n}$ denotes the set of  tuples $\bm{\kappa}=(\kappa_{\ell,j})$ whose components satisfy
$\sum_{\ell=1}^{i}\sum_{j=i}^{n}\kappa_{\ell,j}\leq r-1$ for $1\le i\le n$. 
Then
\begin{align}
&\big\{\Delta\xi^{\left(r\right)}_{\lambda_1}\leq y_1,\ldots, \Delta\xi^{\left(r\right)}_{\lambda_n}\leq y_n\big\}
= \bigcap_{i=1}^n\big\{\Delta\xi^{\left(r\right)}_{\lambda_i}\leq y_i\big\}=  \notag\\
&
\bigcap_{i=1}^n
\bigcup_{\bm{\kappa}=(\kappa_{\ell,j})\in B_{r,n,i}} 
\bigcap_{1\leq \ell\leq j\leq n}\big\{V_{\ell,j}=\kappa_{\ell,j}\big\}
= 
\bigcup_{\bm{\kappa}=(\kappa_{\ell,j})\in A_{r,n}} \bigcap_{1\leq \ell\leq j\leq n}\left\{V_{\ell,j}=\kappa_{\ell,j}\right\}.
\label{eq: set Arn = event xi}
\end{align}
The same  construction  holds with
$\widetilde{V}_{\ell,j,t}$ in place of $V_{\ell,j}$, which means
we can relate $\left\{Z_{r,t,\lambda_i}\leq y_i\right\}$ to $\{\sum_{\ell=1}^{i}\sum_{j=i}^{n}\widetilde{V}_{\ell,j,t}\leq r-1\}$ using the same sets $B_{r,n,i}$.
Thus
\be\label{eq: set Arn = event Z}
\left\{Z_{r,t,\lambda_1}\leq y_1,\ldots, Z_{r,t,\lambda_n}\leq y_n\right\}
=
 \bigcup_{\bm{\kappa}=(\kappa_{\ell,j})\in A_{r,n}}\bigcap_{1\leq \ell\leq j\leq n}\left\{\widetilde{V}_{\ell,j,t}=\kappa_{\ell,j}\right\}
\ee
for the same sets $A_{r,n}$. 

Due to the Poisson nature of the jumps of the  processes $Z$ and $\xi$ in \eqref{eq: set Arn = event xi} and \eqref{eq: set Arn = event Z}, counts of the numbers of points falling in disjoint subrectangles are independent;
in particular, the events $\left\{V_{\ell,j}=\kappa_{\ell,j}\right\}$ 
are independent for all pairs $\ell, j$, $1\leq \ell\leq j\leq n$, and the same is true for the events $\{\widetilde{V}_{\ell,j,t}=
\kappa_{\ell,j}\}$.
Furthermore, the events
\ben 
\bigcap_{1\leq \ell\leq j\leq n}\left\{V_{\ell,j}=\kappa_{\ell,j}\right\}
\quad {\rm and}\quad  \bigcap_{1\leq \ell\leq j\leq n}\left\{V_{\ell,j}=
\overline{\kappa}_{\ell,j}\right\}
\een
 are disjoint if  $\kappa_{\ell,j}\neq \overline{\kappa}_{\ell,j}$
 for at least one  tuple $(\ell,j)$, and the same is true for the tilde version also.

Thus, \eqref{eq: set Arn = event xi} and \eqref{eq: set Arn = event Z} imply
\begin{align}
\mathbb{P}\big(\Delta\xi^{\left(r\right)}_{\lambda_1}\leq y_1,\ldots, \Delta\xi^{\left(r\right)}_{\lambda_n}\leq y_n\big)\label{eq: sum prod P Vij}
=
 \sum_{\bm{\kappa}=(\kappa_{\ell,j})\in A_{r,n}}\prod_{1\leq \ell\leq j\leq n}
\mathbb{P}\left\{V_{\ell,j}=\kappa_{\ell,j}\right\}
\end{align}
and 
\begin{align*}
\mathbb{P}\left(Z_{r,t,\lambda_1}\leq y_1,\ldots, Z_{r,t,\lambda_n}\leq y_n\right)
 =
 \sum_{\bm{\kappa}=(\kappa_{\ell,j})\in A_{r,n}}\prod_{1\leq \ell\leq j\leq n}
 \mathbb{P}\{\widetilde{V}_{\ell,j,t}=\kappa_{\ell,j}\}.
\end{align*}
Hence, to prove \eqref{eq: fin dim distr to show}, 
it remains only to show that for all $m\in\mathbb{N}_0$ the probabilities of the elementary events
$\{\widetilde{V}_{\ell,j,t}
=m\}$
converge to  the probabilities of the events
$\left\{V_{\ell,j}=m\right\}$ 
as $t\dto 0$. 
If we 
define $N_I:\mathbb{R}^+\to\mathbb{N}$ by 
\be\label{Ndef}
N_{I}(x)\coloneqq\#\left\{s\in I\colon \Delta X_s>x\right\},
\ee
 where $I$ is any subinterval of $(0,\infty)$,  
and set 
\begin{align*}
\gamma_{j,t}\coloneqq \overline{\Pi}\big(\overline{\Pi}^{\leftarrow}\big(\left(t y_{j}\right)^{-1}\big)\big)-
\overline{\Pi}\big(\overline{\Pi}^{\leftarrow}\big(\left(t y_{j+1}\right)^{-1}\big)\big), 
\end{align*}
then we can write
\bean
\widetilde{V}_{\ell,j,t}
&=&
N_{t\left(\lambda_{\ell-1},\lambda_{\ell}\right]}
(\pibar\big(\overline{\Pi}^{\leftarrow}\left(t y_{j+1})^{-1}\right)\big)
-
N_{t\left(\lambda_{\ell-1},\lambda_{\ell}\right]}
(\pibar\big(\overline{\Pi}^{\leftarrow}\left(t y_{j})^{-1}\right)\big)   \cr
&&\cr
&\sim&
{\rm Poiss}\big(t\left(\lambda_{\ell}-\lambda_{\ell-1}\right)
\gamma_{j,t}\big).
\eean
Noting further  that
\ben 
 \lim_{t\dto 0}t\gamma_{j,t}=\frac{1}{y_j}-\frac{1}{y_{j+1}},
\een
which follows easily from the slow variation
of $\pibar(x)$ at 0 and the relation
$\overline{\Pi}\left(\overline{\Pi}^{\leftarrow}\left(x\right)\right)
 \leq x <\overline{\Pi}\left(\overline{\Pi}^{\leftarrow}\left(x-\right)
 \right)$,  $x>0$, 
 the convergence of the probabilities of the elementary events finally follows from
\begin{align*}
\MoveEqLeft\lim_{t\dto 0}\mathbb{P}\big(\widetilde{V}_{\ell,j,t}
=m\big)\notag\\
 &=
 \lim_{t\dto 0}
e^{-t\left(\lambda_{\ell}-\lambda_{\ell-1}\right) 
\gamma_{j,t}}
\cdot
\frac{\left(\lambda_{\ell}-\lambda_{\ell-1}\right)^{m}}{m!}\cdot (t\gamma_{j,t})^{m}\notag\\
&=
e^{
-\left(\lambda_{\ell}-\lambda_{\ell-1}\right)
\left(1/y_j-1/y_{j+1}\right)}
\cdot
\frac{\left(\lambda_{\ell}-\lambda_{\ell-1}\right)^{m}}{m!}\cdot\left(\frac{1}{y_{j}}-\frac{1}{y_{j+1}}\right)^{m}\notag\\
&=\mathbb{P}\left(V_{\ell,j}=m\right),
\end{align*}
for all pairs $\ell, j$ fulfilling $1\leq \ell\leq j\leq n$. 
With this, we have completed the proof of  finite dimensional convergence in
Proposition \ref{prop4.3}.\halmos

\subsection{Proof of tightness  in Proposition \ref{prop4.3}}
Recall the $Z_{r,t,\lambda}$ defined in \eqref{Zdef},
which are  positive and  nondecreasing in $\lambda$ for each $r\in\N$ and $t>0$, 
and have the convergence behaviour described in Proposition \ref{prop4.3}. 
In this subsection we show:
\begin{Prop}\label{prop4.4}
Assume $\Pi$ has tail $\overline{\Pi}$  slowly varying at zero. 
Then for all $r\in\mathbb{N}$ the process 
$\big(\big(t\overline{\Pi}\big(\Delta X_{t\lambda}^{(r)}\big)\big)^{-1}\big)_{0<\lambda\leq 1}$ is tight in  $\D[0,1]$ as $t\dto 0$.
\end{Prop}

\medskip\noindent{\bf Proof of Proposition \ref{prop4.4}:}\
We use Theorem 15.3 of Billingsley  (1968), where the result is only stated for discrete time but can immediately be generalised to continuous time as  in the next lemma.

\begin{lemma}\label{4.4}
For each $r\in\N$ the process $(Z_{r,t,\lambda})_{0<\lambda\le 1}$ indexed by $t>0$  is tight
in 
$\D[0,1]$
as $t\dto 0$  if and only if the following conditions hold:

\smallskip
\noindent
(i)\ 
\vskip-0.8cm
 \begin{align}\label{eq: cond 1 billingsley}
 \lim_{y\to\infty}\limsup_{t\dto 0}
  \mathbb{P}\Big(\sup_{0<\lambda\leq 1}Z_{r,t,\lambda}>y\Big)=0;
 \end{align}
(ii)\  for all $y>0$, 
\begin{align}\label{eq: cond 2a billingsley}
 \lim_{\delta\dto 0}\limsup_{t\dto 0}
 \mathbb{P}\Big(
 \sup_{\lambda_1,\lambda_2\in A_{\delta}}
  \sup_{\lambda\in[\lambda_1,\lambda_2]}
\min\left\{Z_{r,t,\lambda}-Z_{r,t,\lambda_1},
Z_{r,t,\lambda_2}-Z_{r,t,\lambda}\right\}>y\Big)\!=\!0,
\end{align}
where 
\begin{align}\label{Adef}
 A_{\delta}:=
 \left\{\lambda_1,\lambda_2 \in (0,1): \, 
 \lambda_1\leq \lambda_2,\, \lambda_2-\lambda_1\leq \delta\right\};
\end{align}
(iii)\  for all $y>0$, 
\begin{align}\label{eq: cond 2b billingsley}
 \lim_{\delta\dto 0}\limsup_{t\dto 0}
 \mathbb{P}\Big(\sup_{\lambda_1,\lambda_2\in\left[0,\delta\right)}
\left|Z_{r,t,\lambda_2}-Z_{r,t,\lambda_1}\right|>y\Big)=0;
\end{align}
(iv)\  for all $y>0$, 
\begin{align}\label{eq: cond 2c billingsley}
 \lim_{\delta\dto 0}\limsup_{t\dto 0}
 \mathbb{P}\Big( \sup_{\lambda_1,\lambda_2\in\left[1-\delta,1\right)}
\left|Z_{r,t,\lambda_2}-Z_{r,t,\lambda_1}\right|>y\Big)=0.
\end{align}
\end{lemma}

In what follows we prove 
\eqref{eq: cond 1 billingsley},
\eqref{eq: cond 2a billingsley}, 
\eqref{eq: cond 2b billingsley} and
\eqref{eq: cond 2c billingsley}
in sequence, keeping $r\in\N$ fixed. 

\medskip\noindent{\bf Proof of Condition (i):}\ 
The probability in the lefthand  side of \eqref{eq: cond 1 billingsley} is
\bea\label{eq: 1-exp}
 \mathbb{P}\big(\sup_{0<\lambda\leq 1} Z_{r,t,\lambda}>y\big)
  &=&
 \mathbb{P}\big(\Delta X_{t}^{(r)}>\overline{\Pi}^{\leftarrow}((ty)^{-1})\big)  \cr
&\le&
 \mathbb{P}\big(\Delta X_{t}^{(1)}>\overline{\Pi}^{\leftarrow}((ty)^{-1})\big)   \cr
 &=&
  1-\mathbb{P}\big(N_{[0,t)}\left(\overline{\Pi}^{\leftarrow}((ty)^{-1})\right)=0\big)   \cr
 &=&
 1-\exp\big(-t\overline{\Pi}\left(\overline{\Pi}^{\leftarrow}((ty)^{-1})\right)\big)   \cr
 &\le& 
 1-\exp\left(-1/y\right).
\eea
(Recall the definition of $N_I$ in \eqref{Ndef}). 
The last inequality in \eqref{eq: 1-exp} follows from the fact that $\overline{\Pi}\left(\overline{\Pi}^{\leftarrow}\left(x\right)\right)\le x$, $x>0$.
Letting $y\to \infty$ in \eqref{eq: 1-exp} gives \eqref{eq: cond 1 billingsley}.

\medskip\noindent{\bf Proof of  Condition   (ii):}\ 
In the following,  keep $y>0$ and $\eta>0$ fixed, and take
$\lambda_0\in(0,1)$ such that $1-e^{-2\lambda_0/y}<\eta/2$.
Recall $A_\delta$ in \eqref{Adef} and define
\ben
 A_{\delta}^\le(\lambda_0) :=
 \left\{\lambda_1,\lambda_2 \in A_\delta:  
 \lambda_1\leq \lambda_0\right\}
 \ {\rm and}\
  A_{\delta}^>(\lambda_0) :=
 \left\{\lambda_1,\lambda_2 \in A_\delta:  
 \lambda_1> \lambda_0\right\}.
\een
Decompose the probability in the lefthand  side of \eqref{eq: cond 2a billingsley}  as
\begin{align}\label{eq: Zrtlambda, 2nd summand}
 &\mathbb{P}\Big(
 \sup_{\lambda_1,\lambda_2\in A_{\delta}}
  \sup_{\lambda\in[\lambda_1,\lambda_2]}
\min\left\{Z_{r,t,\lambda}-Z_{r,t,\lambda_1},
Z_{r,t,\lambda_2}-Z_{r,t,\lambda}\right\}>y\Big)  \nonumber\\
&\leq 
 \mathbb{P}\Big(
 \sup_{\lambda_1,\lambda_2\in A_{\delta}^\le(\lambda_0)}
  \sup_{\lambda\in[\lambda_1,\lambda_2]}
\min\left\{Z_{r,t,\lambda}-Z_{r,t,\lambda_1},
Z_{r,t,\lambda_2}-Z_{r,t,\lambda}\right\}>y\Big)  \nonumber\\
&\ \ \ 
+
\mathbb{P}\Big(
 \sup_{\lambda_1,\lambda_2\in A_{\delta}^>(\lambda_0)}
  \sup_{\lambda\in[\lambda_1,\lambda_2]}
\min\left\{Z_{r,t,\lambda}-Z_{r,t,\lambda_1},
Z_{r,t,\lambda_2}-Z_{r,t,\lambda}\right\}>y\Big).
\end{align}

In the first summand on the  righthand side  of 
\eqref{eq: Zrtlambda, 2nd summand},
$\lambda_1\le \lambda\le \lambda_2\le \lambda_1+\delta \le \lambda_0+\delta$, 
so the probability is bounded above by
\begin{align}\label{eq: lambda2 < lambda0}
\mathbb{P}\big(\sup_{\lambda\in[\lambda_1,\lambda_0+\delta]} Z_{r,t,\lambda}>y\big)
\leq   1-e^{-(\lambda_0+\delta)/y}, 
\end{align}
just as in \eqref{eq: 1-exp}.
When $\delta$ is chosen less than $\lambda_0$, the  righthand side  is 
less than $  1-e^{-2\lambda_0/y}<\eta/2$.

Next  we estimate the second summand on the  righthand side  of \eqref{eq: Zrtlambda, 2nd summand}.  
In it, $\lambda_1>\lambda_0$. Take $\lambda\in[\lambda_1,\lambda_2]$ and $\gamma_1,\gamma_2>0$, and set 
\begin{align*}
 \Gamma_t\coloneqq 
\left\{\gamma_1\leq  Z_{r,t,\lambda_0} 
\le  Z_{r,t,\lambda} \leq \gamma_2\right\}.
\end{align*}
Now \eqref{21} implies 
$  Z_{r,t,\lambda} \todr  \Delta\xi^{(r)}_\lambda$ for each $\lambda\in (0,1]$ 
as $t\dto 0$.
The Cauchy process  $(\xi_\lambda)_{\lambda\ge 0}$ has L\'evy  measure $x^{-2}\rmd x\bfeins_{\{x>0\}}$,
 so the number of jumps exceeding $x>0$ up till time 
$\lambda$ is Poisson with expectation $\lambda/x$.  Thus
\bean 
\PP(\Delta\xi^{(r)}_\lambda\le x) &=&
\PP(\# \{s\in(0,  \lambda):  \Delta \xi_s> x\}\ \le r-1) \cr
&=&
   e^{-\lambda/x}\sum_{j=0}^{r-1}  \frac {(\lambda/x)^j} { j!}.
\eean
This defines a proper distribution with no mass at 0: $\PP(\Delta\xi^{(r)}_\lambda=0)=0$, which is
 continuous as $x\dto 0$.
Thus   we can choose $\gamma_1>0$ small enough, $\gamma_2>0$  large enough and $t_0$ small enough 
so that, for all $t\in(0,t_0)$,
\begin{align}\label{eq: > gamma eta/3} 
\mathbb{P}(\Gamma_t^c)=
\mathbb{P}\left( Z_{r,t,\lambda_0} < \gamma_1\right)
+ \mathbb{P}\left( Z_{r,t,\lambda} >\gamma_2\right)
<\eta/2.
\end{align}

The next task is to show that, for any $t>0$ and $\kappa,\mu$ with $\lambda_0<\kappa<\mu\le 1$, 
\begin{align}\label{Zin}
 \left\{Z_{r,t,\mu}-Z_{r,t,\kappa}>y\right\}
 \cap \Gamma_t
&=
 \left\{\frac{1}{\overline{\Pi}\big(\Delta X^{\left(r\right)}_{t\mu}\big)}-\frac{1}{\overline{\Pi}\big(\Delta X^{\left(r\right)}_{t\kappa}\big)}\geq t y\right\}\cap \Gamma_t\nonumber \\
&\subseteq
\left\{\frac{\Delta X^{(r)}_{t\mu}-\Delta X^{\left(r\right)}_{t\kappa}}{\Delta X^{\left(r\right)}_{t\lambda_0}\ 
\overline{\Pi}\big(\Delta X^{\left(r\right)}_{t}\big)}\geq t y\right\}\cap \Gamma_t.
\end{align}
We again apply   Potter's bounds, see  Theorem 1.5.6 of Bingham, Goldie and Teugels (1987),
 where the theorem is stated for functions slowly varying at infinity but can be immediately transferred to functions slowly varying at zero. 
In one of its forms it states that for a function $L$ slowly varying at  zero there exists $T>0$ such that
\begin{align}\label{uv}
\min\left\{\frac{u}{v},\frac{v}{u}\right\}< \frac{L\left(u\right)}{L\left(v\right)}<\max\left\{\frac{u}{v},\frac{v}{u}\right\},
\end{align}
for all $u,v\in (0,T]$.
On $\Gamma_t$ we have $ Z_{r,t,\lambda_0} =1/\big(t \pibar\big( \Delta X^{(r)}_{t\lambda_0}\big)\big)\le \gamma_1$, so
\be\label{40a}
\pibarinv\Big(\frac{1}{\gamma_1t}\Big)
\le \Delta X^{(r)}_{t\lambda_0}
\le \Delta X^{(r)}_{t\kappa}.
\ee
Choosing $0<t\le 1/(\gamma_1\pibar(T))$, we have 
$\pibarinv(1/(\gamma_1t))\le T$, so by \eqref{uv}, 
\ben
\frac{\pibar\big( \Delta X^{(r)}_{t\mu}\big)}
{\pibar\big( \Delta X^{(r)}_{t\kappa}\big)}
\ge \frac{\Delta X^{(r)}_{t\kappa}}{\Delta X^{(r)}_{t\mu}}.
\een
This yields
\begin{align*}
\frac{1}{\overline{\Pi}\big(\Delta X^{\left(r\right)}_{t\mu}\big)}-\frac{1}{\overline{\Pi}\big(\Delta X^{\left(r\right)}_{t\kappa}\big)}
&\leq 
\frac{1}{\overline{\Pi}\big(\Delta X^{\left(r\right)}_{t\mu}\big)}
-
\left(\frac{\Delta X^{\left(r\right)}_{t\kappa}}{\Delta X^{\left(r\right)}_{t\mu}}\right)
\frac{1}{\overline{\Pi}\big(\Delta X^{\left(r\right)}_{t\mu}\big)}\notag\\
&=
\frac{\Delta X^{\left(r\right)}_{t\mu}-\Delta X^{\left(r\right)}_{t\kappa}}
{\Delta X^{\left(r\right)}_{t\mu}\ \overline{\Pi}\big(\Delta X^{(r)}_{t\mu}\big)}\notag\\
&\leq 
\frac{\Delta X^{\left(r\right)}_{t\mu}-\Delta X^{\left(r\right)}_{t\kappa}}
{\Delta X^{\left(r\right)}_{t\lambda_0}\ \overline{\Pi}\big(\Delta X^{\left(r\right)}_{t}\big)},
\end{align*}
on the event $\Gamma_t$, when $0<t\le 1/(\gamma_1\pibar(T))$.
With this inequality we have proved the inclusion in \eqref{Zin}.

Continuing from \eqref{Zin},  argue from \eqref{40a} that, on $\Gamma_t$, 
\begin{align}\label{eq: Delta Xtmu - Delta Xtkappa}
\MoveEqLeft\left\{\frac{\Delta X^{\left(r\right)}_{t\mu}-\Delta X^{\left(r\right)}_{t\kappa}}
{\Delta X^{\left(r\right)}_{t\lambda_0}\ \overline{\Pi}\big(\Delta X^{(r)}_{t}\big)}>t y\right\}    
\subseteq
\left\{\Delta X^{\left(r\right)}_{t\mu}-\Delta X^{\left(r\right)}_{t\kappa}>\frac{y}{\gamma_2} \overline{\Pi}^{\leftarrow}\left(\frac{1}{t\gamma_1}\right)\right\}
\end{align}
(here note too that 
$1/ \overline{\Pi}\big(\Delta X^{(r)}_{t}\big)
= tZ_{r,t,1}\le t \gamma_2$ on $\Gamma_t$).
For the following, set $a_t\coloneqq (y/\gamma_2) \overline{\Pi}^{\leftarrow}\left(1/\left(t\gamma_1\right)\right)$.
Applying  \eqref{eq: Delta Xtmu - Delta Xtkappa} once for $\mu\coloneqq \lambda$ and $\kappa\coloneqq \lambda_1$,
 and once for $\mu\coloneqq \lambda_2$ and $\kappa\coloneqq \lambda$, we obtain 
 \begin{align*}
 &\mathbb{P}\Big(
 \sup_{\lambda_1,\lambda_2\in A_{\delta}^>(\lambda_0)}
  \sup_{\lambda\in[\lambda_1,\lambda_2]}
\min\left\{Z_{r,t,\lambda}-Z_{r,t,\lambda_1},
Z_{r,t,\lambda_2}-Z_{r,t,\lambda}\right\}>y;\, \Gamma_t\Big) 
\nonumber\\
&\leq 
 \mathbb{P}\Big(
 \sup_{\lambda_1,\lambda_2\in A_{\delta}^>(\lambda_0)}
  \sup_{\lambda\in[\lambda_1,\lambda_2]}
\min\big(\Delta X^{\left(r\right)}_{t\lambda}-\Delta X^{\left(r\right)}_{t\lambda_1},\Delta X^{\left(r\right)}_{t\lambda_2}-\Delta X^{\left(r\right)}_{t\lambda}\big)>a_t;\, \Gamma_t\Big).
\end{align*}
For given $\lambda_1,\lambda_2$, the event
\begin{align}\label{Ev}
\Big\{
 \sup_{\lambda\in\left[\lambda_1, \lambda_2\right]}
 \min\big(\Delta X^{\left(r\right)}_{t\lambda}-\Delta X^{\left(r\right)}_{t\lambda_1},\Delta X^{\left(r\right)}_{t\lambda_2}-\Delta X^{\left(r\right)}_{t\lambda}\big)>a_t\Big\}
\end{align}
requires that there exist at least two points $s_1,s_2\in \left(\lambda_1, \lambda_2\right]$ such that 
$\Delta X_{ts_1}>a_t$ and
$\Delta X_{ts_2}>a_t$. 
To see this, assume there is no point $s\in\left(\lambda_1,\lambda_2\right]$ with $\Delta X_{ts}>a_t$.
Then $\Delta X^{\left(r\right)}_{t\lambda_2}-\Delta X^{\left(r\right)}_{t\lambda_1}\le a_t$ and thus
$\Delta X^{\left(r\right)}_{t\lambda}-\Delta X^{\left(r\right)}_{t\lambda_1}\le a_t$ and $\Delta X^{\left(r\right)}_{t\lambda_2}-\Delta X^{\left(r\right)}_{t\lambda}\le a_t$ hold for any $\lambda\in \left(\lambda_1,\lambda_2\right]$. 
This is not possible under \eqref{Ev}.
If there is only one point $s\in\left(\lambda_1,\lambda_2\right]$ with $\Delta X_{ts}>a_t$, then for any $\lambda\in \left(\lambda_1,\lambda_2\right]$ we have that either
$\Delta X^{\left(r\right)}_{t\lambda}-\Delta X^{\left(r\right)}_{t\lambda_1}\le a_t$ or $\Delta X^{\left(r\right)}_{t\lambda_2}-\Delta X^{\left(r\right)}_{t\lambda}\le a_t$, also  not possible under \eqref{Ev}. 
Hence we deduce
\begin{align}\label{eq: Pt x sup}
&\mathbb{P}\Big(
 \sup_{\lambda_1,\lambda_2\in A_{\delta}^>(\lambda_0)}
  \sup_{\lambda\in[\lambda_1,\lambda_2]}
\min\left\{Z_{r,t,\lambda}-Z_{r,t,\lambda_1},
Z_{r,t,\lambda_2}-Z_{r,t,\lambda}\right\}>y;\, \Gamma_t\Big) 
\nonumber\\
&\leq 
\mathbb{P}\left(\exists\lambda\in\left[\lambda_0, 1-\delta\right]\colon N_{\left[t\lambda,t\left(\lambda+\delta\right)\right)}\left(a_t\right)\geq 2; \, \Gamma_t\right).
\end{align}

Now   define intervals $I_{k,t,\delta}\coloneqq \left[t(\lambda_0+k\delta),\,
t(\lambda_0+(k+2)\delta)\right)$, for $t>0$, $\delta>0$ and $k\in\N$. Note that the length of each of these intervals is $2t\delta$.
Further,  define the integers 
$k_{\delta}\coloneqq\left\lceil \left(1-\lambda_0\right)/\delta\right\rceil$.

For given $\delta>0$ and $\lambda\in\left[\lambda_0, 1-\delta\right]$ there exists $k\in\left[0,k_{\delta}\right]\cap\N$ 
such that $\lambda\in\left[\lambda_0+k\delta, \lambda_0+\left(k+1\right)\delta\right)$,
hence
$t \lambda\in\left[t\left(\lambda_0+k\delta\right), t\left(\lambda_0+\left(k+1\right)\delta\right)\right)$.
This implies 
$\left[t \lambda, t\left(\lambda+\delta\right)\right)
\subseteq\left[t\left(\lambda_0+k\delta\right), t\left(\lambda_0+\left(k+1\right)\delta\right)\right)=I_{k,t,\delta}$
for the same $k$, so for each interval $\left[t \lambda, t\left(\lambda+\delta\right)\right)$ there exists 
$k\in\left[0,k_{\delta}\right]\cap\N$ such that $\left[t \lambda, t\left(\lambda+\delta\right)\right)\subset I_{k,t,\delta}$.

Thus,
\be\label{eq: exists lambda}
\left\{\exists \lambda\in\left[\lambda_0, 1-\delta\right]\colon N_{\left[t\lambda,t\left(\lambda+\delta\right)\right)}\left(a_t\right)\geq 2\right\}
\subseteq \bigcup_{k=0}^{k_{\delta}} \left\{ N_{I_{k,t,\delta}}\left(a_t\right)\geq 2\right\}.
\ee
The intervals $I_{k,t,\delta}$  are constructed in such a way that every second interval is disjoint from the preceding one.
Thus the events 
$\left\{ N_{I_{2k-1,t,\delta}}\left(a_t\right)>2\right\}$ for $k\in\left[1,\left\lceil k_{\delta}/2\right\rceil\right]\cap \N$
are mutually independent, as are the events $\left\{ N_{I_{2k,t,\delta}}\left(a_t\right)>2\right\}$ for
 $k\in\left[0,\left\lfloor k_{\delta}/2\right\rfloor\right]\cap \N_0$.
Accordingly, write the  righthand side  of \eqref{eq: exists lambda} as
\begin{align*}
\bigcup_{k=0}^{k_{\delta}} \left\{ N_{I_{k,t,\delta}}\left(a_t\right)\geq 2\right\}
&\!=\bigcup_{k=0}^{\left\lfloor k_{\delta}/2\right\rfloor} \left\{ N_{I_{2k,t,\delta}}\left(a_t\right)\geq 2\right\}
\cup \bigcup_{k=1}^{\left\lceil k_{\delta}/2\right\rceil} \left\{ N_{I_{2k-1,t,\delta}}\left(a_t\right)\geq 2\right\},
\end{align*}
and combine this with \eqref{eq: Pt x sup} and \eqref{eq: exists lambda} to get 
\begin{align*}
\MoveEqLeft\mathbb{P}\Big(
 \sup_{\lambda_1,\lambda_2\in A_{\delta}^>(\lambda_0)}
  \sup_{\lambda\in[\lambda_1,\lambda_2]}
\min\left\{Z_{r,t,\lambda}-Z_{r,t,\lambda_1},\, Z_{r,t,\lambda_2}-Z_{r,t,\lambda} \right\}\geq y;\,
 \Gamma_t\Big)\\ 
&\leq \mathbb{P}\bigg(\bigcup_{k=0}^{\left\lfloor k_{\delta}/2\right\rfloor} \left\{ N_{I_{2k,t,\delta}}\left(a_t\right)\geq 2\right\}\bigg)
+\mathbb{P}\bigg(\bigcup_{k=1}^{\left\lceil k_{\delta}/2\right\rceil} \left\{ N_{I_{2k-1,t,\delta}}\left(a_t\right)\geq 2\right\}\bigg)\\
&= 2- \mathbb{P}\bigg(\bigcup_{k=0}^{\left\lfloor k_{\delta}/2\right\rfloor} \left\{ N_{I_{2k,t,\delta}}\left(a_t\right)\leq 1\right\}\bigg)
-\mathbb{P}\bigg(\bigcup_{k=1}^{\left\lceil k_{\delta}/2\right\rceil} \left\{ N_{I_{2k-1,t,\delta}}\left(a_t\right)\leq 1\right\}\bigg).
\end{align*}
The events $\left\{ N_{I_{2k,t,\delta}}\left(a_t\right)\leq 1\right\}$ with  $k\in\left[0,\left\lfloor k_{\delta}/2\right\rfloor\right]\cap \N_0$ 
are mutually independent as well as the events $\left\{ N_{I_{2k-1,t,\delta}}\left(a_t\right)\leq 1\right\}$ with  $k\in\left[0,\left\lceil k_{\delta}/2\right\rceil\right]\cap \N$.
This implies
\begin{align}
\MoveEqLeft\mathbb{P}\Big(
 \sup_{\lambda_1,\lambda_2\in A_{\delta}^>(\lambda_0)}
  \sup_{\lambda\in[\lambda_1,\lambda_2]}
\min\left\{Z_{r,t,\lambda}-Z_{r,t,\lambda_1},\, Z_{r,t,\lambda_2}-Z_{r,t,\lambda} \right\}\geq y;\,
 \Gamma_t\Big)\notag\\ 
&\leq  2-\bigg(\prod_{k=0}^{\left\lfloor k_{\delta}/2\right\rfloor}
\mathbb{P}\left( N_{I_{2k,t,\delta}}\left(a_t\right)\leq 1\right)\bigg)
-\bigg(\prod_{k=1}^{\left\lceil k_{\delta}/2\right\rceil}
\mathbb{P}\left( N_{I_{2k-1,t,\delta}}\left(a_t\right)\leq 1\right)\bigg)\notag\\
&= 2\left(1-\mathbb{P}^{\left\lfloor k_{\delta}/2\right\rfloor+1}\left( N_{I_{0,t,\delta}}\left(a_t\right)\leq 1\right)\right).\label{eq: 1-P Nk}
\end{align}
Here the last equality follows from the fact that each of the intervals has the same length
and thus the probabilities $\mathbb{P}\left( N_{I_{k,t,\delta}}\left(a_t\right)\leq 1\right)$
are equal for all $k\in\left[0,k_{\delta}\right]\cap\N_0$.
Furthermore
\ben
\mathbb{P}\left( N_{I_{0,t,\delta}}\left(a_t\right)\leq 1\right)
=
\left(1+2t\delta \pibar(a_t))\right)
e^{-2t\delta \overline{\Pi}(a_t)},
\een
thus,
\be\label{kd}
\mathbb{P}^{\left\lfloor k_{\delta}/2\right\rfloor+1}  \left( N_{I_{0,t,\delta}}(a_t)\leq 1\right)
=(1+ \delta c_t)^{\left\lfloor k_{\delta}/2\right\rfloor+1}  e^{-\delta(\left\lfloor k_{\delta}/2\right\rfloor+1) c_t},
\ee
where $c_t:= 2t\pibar(a_t)$.
Letting $t\dto 0$, so that
\ben
c_t=2t\pibar(a_t)=2t\pibar(t\pibar^{\leftarrow}(1/(t \gamma_1))y/\gamma_2))\to2/\gamma_1,
\een
 followed by  $\delta\dto 0$, so that
$(\left\lfloor k_{\delta}/2\right\rfloor+1)\delta\to (1-\lambda_0)/2$, shows that the  righthand side  of \eqref{kd} tends to 
$e^{(1-\lambda_0)/\gamma_1}e^{-(1-\lambda_0)/\gamma_1}=1$.
Then we deduce from \eqref{eq: 1-P Nk} that
\ben 
\limsup_{t\dto 0}
\mathbb{P}\Big(
 \sup_{\lambda_1,\lambda_2\in A_{\delta}^>(\lambda_0)}
  \sup_{\lambda\in[\lambda_1,\lambda_2]}
\min\left\{Z_{r,t,\lambda}-Z_{r,t,\lambda_1},\, Z_{r,t,\lambda_2}-Z_{r,t,\lambda} \right\}>y;\,
 \Gamma_t\Big) 
\een
tends to 0 as $\delta\dto 0$.
Combining this with \eqref{eq: lambda2 < lambda0} and 
 \eqref{eq: > gamma eta/3}  yields that
\ben
\lim_{\delta\dto 0}\limsup_{t\dto 0}
\mathbb{P}\Big(
 \sup_{\lambda_1,\lambda_2\in A_{\delta}^>(\lambda_0)}
  \sup_{\lambda\in[\lambda_1,\lambda_2]}
\min\left\{Z_{r,t,\lambda}-Z_{r,t,\lambda_1},\, Z_{r,t,\lambda_2}-Z_{r,t,\lambda} \right\}>y\Big)  
\een
is less than $\eta$, 
and since $\eta$ is arbitrary this proves  \eqref{eq: cond 2a billingsley}. \halmos

\medskip\noindent{\bf Proof of   Condition  (iii):}\ 
The probability in the lefthand  side of  \eqref{eq: cond 2b billingsley} can be written as 
\ben
\mathbb{P}\Big(\sup_{\lambda_1,\lambda_2\in\left[0,\delta\right)}\left|Z_{r,t,\lambda_2}-Z_{r,t,\lambda_1}\right|>y\Big),
\een
 and this is no larger than 
$\mathbb{P}\big(Z_{r,t,\delta/2}>y\big)$. 
Using a similar calculation  as in \eqref{eq: 1-exp},  there exists
$t_5>0$ such that for all $t\in(0,t_5)$ 
\begin{align*}
 \mathbb{P}\big(Z_{r,t,\delta}>y\big) \leq 1-e^{-\delta/y}.
\end{align*}
For  fixed $y>0$ and  $\delta>0$ small enough this is no larger than $\eta$.

\medskip\noindent{\bf Proof of   Condition  (iv):}\ 
The probability in the lefthand  side of \eqref{eq: cond 2c billingsley} 
is no larger than 
$ \mathbb{P}\big(Z_{r,t,1}-Z_{r,t,1-\delta}>y\big)$. 
Just as in the proof of Condition (ii) there exists $t_6>0$ such that for $t\in(0,t_6)$  
\begin{align*}  
\mathbb{P}\Big(
Z_{r,t,1}-Z_{r,t,1-\delta}>y\Big)
 &\leq \mathbb{P}\left(N_{\left[t\left(1-\delta\right),t\right]}\left(a_t\right)\geq 1\right)
 +\mathbb{P}\left(\Gamma_t^c\right) \cr
  &\le
1-e^{-\delta t\pibar(a_t)}
  +\eta \cr
 &\le
1-e^{-2\delta/\gamma_1}  +\eta.
\end{align*} 
For   $\delta>0$ small enough this is no larger than $2\eta$,  so the  proof of   Condition  (iv) is complete,
 and this finally completes the proof of Proposition \ref{prop4.4}.
\halmos

\medskip\noindent{\bf Remarks.}\ 
The almost sure condition \eqref{eq: Xt Delta Xt as} may seem anomalous in the midst of the other weak convergence conditions, 
but it is not excessive in context. It follows from the proof of
 Lemma \ref{lem2}, via  \eqref{eq: t Pi X - t Pi Delta X 2a} 
that, with 
 \ben
 R_{t\lambda}:= 
 \frac{  \pibar(\Delta X_{t\lambda}^{(r+1)}) }{\pibar(\upperl{r}{}{X}{t\lambda}{})} \ge  1,
 \een
we have
\ben
\lim_{t\dto 0}
\mathbb{P}\big(\sup_{ 0< \lambda\le 1} R_{t\lambda}>1+\veps\big)=0,
\een
for all $\veps>0$.
Hence, for given $\veps>0$, $\delta>0$,  and $t$ small enough, $t\le t_0(\veps,\delta)$,
 \ben 
\mathbb{P}\big( R_{t\lambda}>1+\veps\ {\rm for\ some}\ \lambda\in(0,1]\big) \le \delta.
\een
But this implies 
 \ben 
\mathbb{P}\big( R_s>1+\veps\ {\rm for\ some}\ s\le t\big) \le \delta
\een
whenever $t\le t_0$.  Hence
\ben
\lim_{t\dto 0} 
 \frac{  \pibar(\Delta X_{t\lambda}^{(r+1)}) }{\pibar(\upperl{r}{}{X}{t\lambda}{})} =1,\ {\rm a.s.}
 \een
This is close to \eqref{eq: Xt Delta Xt as} but  does not imply it in general because the converse  part of Lemma \ref{lem1}
is not true for $\alpha=0$ in general (take, for example, $\pibar(x)=|\log x|$, $f_t=t|\log t|$, $g_t=t$, for $0<x,t<1$).
So we have to impose \eqref{eq: Xt Delta Xt as}  as a side condition.

We remark incidentally that  the slow variation of $\pibar(x)$ at 0
is equivalent to  a weak version of \eqref{eq: Xt Delta Xt as}, namely that 
$\upperl{r}{}{X}{t}{}/ \Delta X_{t}^{(r+1)}\topr 1$
as $t\dto 0$ for $r\in\N_0$  (see Buchmann, Ipsen, Maller (2016). 
A  necessary and sufficient condition for  \eqref{eq: Xt Delta Xt as}  itself in the case $r=1$ is in Maller (2016). 

\section{Convergence of the Trimmed Stable as $\alpha\dto 0$ }\label{s5}
In this section, to complete Figure \ref{diagram} we prove that $\big({}^{(r)}S_\alpha(\lambda)\big)^{\alpha}$
converges to $(\Delta\xi_{\lambda}^{\left(r+1\right)})$ in $\D[0,1]$ as $\alpha\dto 0$, for each $r\in\N_0$. 
First suppose $r=0$.
As in the proof of Kasahara (1986), 
we obtain that $S_\alpha(\lambda)$ can be written as
\begin{align*}
S_\alpha(\lambda)
=\int_{u\in(0, \lambda]}\int_{x>0}x^{1/\alpha} N\left(\mathrm{d}u, \mathrm{d}x\right),
\end{align*}
where $N\left(\mathrm{d}u, \mathrm{d}x\right)$ is a Poisson random measure with intensity measure $\mathrm{d}u\times x^{-2}\mathrm{d}x$.
This is the Poisson random measure governing the jumps of the Cauchy process $\left(\xi_{\lambda}\right)_{0<\lambda\leq 1}$, so we can write
\ben
\big(S_\alpha(\lambda)\big)_{0<\lambda\le 1} 
=\bigg(\sum_{0<s\leq \lambda}\left(\Delta\xi_{\lambda}\right)^{1/\alpha}\bigg)_{0<\lambda\le 1}.
\een
 The  jumps up till  time $\lambda$ of $\xi_\lambda$ can be ordered as
 $\Delta \xi_\lambda^{(1)}\ge \Delta \xi_\lambda^{(2)}\ge \cdots$, and then 
since raising to the power $1/\alpha$ does not change the order of the jumps,
\be\label{prh}
({}^{(r)}S_\alpha(\lambda))^\alpha
=\bigg(\sum_{i\geq r+1} \big(\Delta\xi_{\lambda}^{\left(i\right)}\big)^{1/\alpha}\bigg)^{\alpha}. 
\ee
Using a  classical argument\footnote{When $a_r\ge a_{r+1}\ge \cdots>0$ and $\sum_{i\ge r} a_i<\infty$, then
\ben
\alpha \log \big(\sum_{i\ge r} a_i^{1/\alpha}\big) 
=\log a_r +\alpha\big(1+\sum_{i> r}(a_i/a_r)^{1/\alpha}\big).
\een
Take $\alpha<1$ and choose $i_0(r)>r$ so that $(a_{i_0}/a_r)^{1/\alpha-1}<1$. Then the second term on the righthand is less than
$\alpha\big(i_0-r +\sum_{i> i_0}a_i/a_r\big)\to 0$ as $\alpha\dto 0$.}
we can show that when $\alpha\dto 0$ each term in the process on the  righthand side 
of \eqref{prh}  converges {\it surely}
(i.e., for each $\omega\in \Omega$) to
\begin{align*}
 \sup_{i\geq r+1} \Delta\xi_{\lambda}^{\left(i\right)}=
 \Delta\xi_{\lambda}^{\left(r+1\right)}.
\end{align*}
Consequently,  also the process on the  righthand side  of \eqref{prh}  converges surely
to the process $(\Delta\xi_{\lambda}^{\left(r+1\right)})$.
This of course also implies convergence in distribution.
So we obtain the required result.
\halmos

\subsection*{Acknowledgements}
This research  was partially supported by ARC Grant DP160103037.


\section{Appendix} \label{s6} 
 To give a concrete formula for the finite dimensional distribution for the $r$th jump of a subordinator $(Y_t)$ without drift is tedious in general. 
 However, in case $r=1$ this is a classical result and can be found for example in Chapter 4.1 of Resnick (2008). 
 Let $\Lambda$ be the L\'evy measure of $Y$, with tail $\lambar$.  Then we have,
 for $\lambda_1<\cdots<\lambda_n$ and $y_1<\cdots<y_n$,
 \begin{align*}
  \MoveEqLeft
  \mathbb{P}\big(
\Delta Y_{\lambda_1}^{\left(1\right)}\leq y_1,\ldots, 
\Delta Y_{\lambda_n}^{\left(1\right)}\leq y_n\big)\\
  &=
  e^{-\lambda_1\overline{\Lambda}\left(y_1\right)} e^{-(\lambda_2-\lambda_1)\overline{\Lambda}\left(y_2\right)}
\cdots
e^{-\left(\lambda_n-\lambda_{n-1}\right)
\overline{\Lambda}\left(y_n\right)}.
 \end{align*}
 In case $(Y_t)=(\xi_t)$ is  a Cauchy process, this simplifies to
 \begin{align}\label{sC}
  \MoveEqLeft
  \mathbb{P}\big(
  \Delta\xi_{\lambda_1}^{\left(1\right)}\leq y_1,\ldots, \Delta\xi_{\lambda_n}^{\left(1\right)}\leq y_n\big)\cr
  &=
e^{-\lambda_1/y_1} e^{-(\lambda_2-\lambda_1)/y_2} \cdots e^{-(\lambda_n-\lambda_{n-1})/y_n}.
 \end{align}

We can also get this as a calculation from  \eqref{eq: sum prod P Vij}.
In case $r=1$ we take $A_{1,n}=\big\{\left(\kappa_{\ell,j}\right)_{1\leq \ell\leq j\leq n}\colon \kappa_{\ell,j}=0\big\}$, and the formula simplifies to
\begin{align*}
 \MoveEqLeft
  \mathbb{P}\big(
  \Delta\xi_{\lambda_1}^{\left(1\right)}\leq y_1,\ldots, \Delta\xi_{\lambda_n}^{\left(n\right)}\leq y_n\big)
  =\prod_{1\leq \ell\leq j\leq n}\mathbb{P}\left(V_{\ell,j}=0\right)   \\ 
  &=
  \prod_{1\leq \ell\leq j\leq n}
e^{-(\lambda_\ell-\lambda_{\ell-1})\left(1/y_j-1/y_{j+1}\right)},
\end{align*}
which is the same as the righthand side of \eqref{sC}
(recall  $\lambda_0=0$ and $y_{n+1}=\infty$).

\subsection{A formula for the fidi distribution  of the 2nd largest jump}
 In the following we derive an  implicit formula for $r=2$. 
For larger $r$ the formula could be derived in a similar way.
Let $Y=(Y_t)$ be any subordinator without drift and with L\'evy measure $\Lambda$. 
We aim to give a formula for
\begin{align*}
 \mathbb{P}\big(\Delta Y_{\lambda_1}^{\left(2\right)}<y_1,\ldots,\Delta Y_{\lambda_n}^{\left(2\right)}<y_n\big),
\end{align*}
where $0<\lambda_1<\cdots<\lambda_n$ and $0<y_1<\cdots <y_n$. 
 
Analogously to \eqref{eq: def V ij}, we will set 
 $\lambda_0\coloneqq 0$ and $y_{n+1}\coloneqq\infty$ and then
\begin{align*}
V_{\ell,j}&\coloneqq\#\left\{s\in\left[\lambda_{\ell-1},\lambda_\ell\right)\colon \Delta Y_s \in\left[y_j,y_{j+1}\right)\right\}.
\end{align*}
One way to calculate the finite dimensional distribution would be to construct the set $A_{r,n}$ given in Section \ref{sec: fin dim conv}.
However, this would require constructing the set of triangular arrays fulfilling $\sum_{\ell=1}^i\sum_{j=i}^n\kappa_{\ell,j}\leq r$ for all $i\in\left\{1,\ldots, n\right\}$ at the same time.
To our knowledge there is no simple way to do that. 
So we choose a slightly different approach.

To start, we set  $D_{n+1,n}=\Omega$ and, for $2\leq i\leq n$,
\begin{align*}
D_{i,n}\coloneqq\big\{\Delta Y_{\lambda_i-\lambda_{i-1}}^{\left(2\right)}<y_i,\ldots,\Delta Y_{\lambda_n-\lambda_{n-1}}^{\left(2\right)}<y_n\big\}.
\end{align*}
Then note that 
\begin{align}
\mathbb{P}\big(\Delta Y_{\lambda_i}^{\left(2\right)}<y_i,\ldots,\Delta Y_{\lambda_n}^{\left(2\right)}<y_n\big| \Delta Y_{\lambda_{i-1}}^{\left(1\right)}<y_i\big)
&=\mathbb{P}\left(D_{i,n}\right).\label{eq: def pin}
\end{align}
This follows from the fact that the numbers of jumps in different intervals are independent.
Given there are no jumps exceeding $y_i$ in the interval $\left[0,\lambda_{i-1}\right)$, 
there are in  particular no jumps exceeding $y_{i},\ldots,y_n$. 
Hence, under the condition $\Delta Y_{\lambda_{i-1}}^{\left(1\right)}<y_i$,
the number of jumps exceeding $y_i,\ldots,y_n$ on the intervals
$[0,\lambda_i],\ldots,[0,\lambda_n]$ is the same as
the number of jumps exceeding $y_i,\ldots,y_n$ on the intervals
$(\lambda_{i-1},\lambda_i],\ldots, (\lambda_{i-1},\lambda_n]$.
Since the increments are stationary, we obtain the formula in \eqref{eq: def pin}.

Next we state our recursive formula and give an explanation  following it.
The formula is 
\begin{align}\label{eq: fin dim distr r=2}
 \MoveEqLeft
 \mathbb{P}\big(\Delta Y_{\lambda_1}^{\left(2\right)}<y_1,\ldots,\Delta Y_{\lambda_n}^{\left(2\right)}<y_n\big)
=\prod_{j=1}^n\mathbb{P}\left(V_{1,j}=0\right)
\cdot \mathbb{P}\left(D_{2,n}\right)+
 \notag\\
 &
 \sum_{i=1}^n\mathbb{P}\left(V_{1,i}=1\right)
 \Big(\prod_{i\neq j} \mathbb{P}(V_{1,j}=0)  \Big)
 \Big(\prod_{\ell=2}^i\, \prod_{j=\ell}^n
 \mathbb{P}(V_{\ell,j}=0)\Big)\cdot \mathbb{P}\left(D_{i+1,n}\right),
\end{align}
where by convention we set $\prod_{k=2}^1=1$.
Note that the formula is recursive in that sense 
that the probability of the elementary events $D_{n,n}$ can immediately be calculated by noticing that $D_{n,n}=\left\{V_{n,n}\leq 1\right\}$
and
\begin{align*}
\mathbb{P}\left(V_{\ell,j}=k\right)
&=
e^{-(\lambda_\ell-\lambda_{\ell-1}) \Lambda\left(\left[y_{j-1},y_{j}\right)\right)}
\cdot \frac{\left(\lambda_\ell-\lambda_{\ell-1}\right)^k\cdot \Lambda\left(\left[y_{j-1},y_{j}\right)\right)^k}{k!}.
\end{align*}
For $i<n$ the events $D_{i,n}$ are of the form 
of the lefthand side of \eqref{eq: fin dim distr r=2} with smaller $n$ which specifies the  recursion. 

Notice also 
that $\Delta Y^{\left(2\right)}_{\lambda_1}<y_1$ if and only if 
$\sum_{j=1}^nV_{1,j}\leq 1$.
First assume $\sum_{j=1}^nV_{1,j}=0$. 
Then it suffices to have $\sum_{\ell=2}^k\sum_{j=k}^nV_{\ell,j}\leq 1$ for all $k\in\left\{2,\ldots,n\right\}$.
This is equivalent to $D_{2,n}$ and gives the first summand of \eqref{eq: fin dim distr r=2}.

To obtain the second summand of \eqref{eq: fin dim distr r=2} let us assume that $\sum_{j=1}^n V_{1,j}=1$, 
which is equivalent to the statement that there exists $i\in\left\{1,\ldots,n\right\}$ such that
$V_{1,i}=1$ and $V_{1,\ell}=0$ for all $\ell\neq i$ which
are represented in the sum in \eqref{eq: fin dim distr r=2}. 
Assume that this is the case and 
remember from Section \ref{sec: fin dim conv} 
that $\{\Delta Y_{\lambda_k}^{(2)}<y_k\}=\{\sum_{\ell=1}^k\sum_{j=k}^n V_{\ell,j}\leq 1\}$.
Then in order that $\bigcap_{k=1}^i\{\Delta Y_{\lambda_k}^{(2)}<y_k\}$ holds it is necessary and sufficient 
that 
$\sum_{\ell=2}^{k}\sum_{j=k}^n V_{\ell, j}=0$ for all $k\in\left\{2,\ldots,i\right\}$.
This in turn is equivalent to $V_{\ell,j}=0$ for all pairs $\ell,j$ with $\ell\in\left\{2,\ldots,i\right\}$ and $j\in\left\{\ell,\ldots, n\right\}$.

Given this is the case, then for each of the events $\{\Delta Y_{\lambda_k}^{(2)}<y_k\}=\{\sum_{\ell=1}^k\sum_{j=k}^n V_{\ell,j}\leq 1\}$ with $k\in\left\{i+1,\ldots,n\right\}$ to hold it is additionally necessary and sufficient that 
$\{\sum_{\ell=i+1}^k\sum_{j=k}^n V_{\ell,j}\leq 1\}$ for all $k\in\left\{i+1,\ldots,n\right\}$.
The intersection over the last events with indices $k\in\left\{i+1,\ldots,n\right\}$ is equivalent to
$D_{i+1,n}$.

Combining all these argumentations gives the formula in \eqref{eq: fin dim distr r=2}.
\halmos


\end{document}